\documentclass[12pt]{article}
\usepackage{amsfonts}
\usepackage{amsmath}
\usepackage{url}
  \newcommand{\lab}[1]{\label{#1}}                

 \def\bigpage {
\addtolength{\topmargin}{-2  cm}
\addtolength{\oddsidemargin}{-2  cm}
\addtolength{\textwidth}{+4 cm}
\addtolength{\textheight}{+4 cm} }

\bigpage

\def\proof{\noindent{\bf Proof}\ \ }
\def\qed{~~\vrule height8pt width4pt depth0pt}

\newcommand{\be}{\begin{equation}}
\newcommand{\ee}{\end{equation}}
\newcommand{\bea}{\begin{eqnarray}}
\newcommand{\eea}{\end{eqnarray}}

\newcommand{\bean}{\begin{eqnarray*}}
\newcommand{\eean}{\end{eqnarray*}}
\newcommand\eqn[1]{(\ref{#1})}
\newcommand{\bel}[1]{\be\lab{#1}}

\catcode`@=11
\@addtoreset{equation}{section}

\catcode`@=12

\newtheorem{thm}{Theorem}[section]

\newtheorem{lemma}[thm]{Lemma}

\def\eps{\epsilon}
\def\G{{\cal G}}
\def\C{{\cal C}}
\def\P{{\cal P}}
\def\dv{{\bf d}}
\def\ex{{\bf E\,}}
\def\pr{{\bf P}}

\ifx\hyperlink\undefined
\newcommand{\hrlb}{ }

\newcommand{\hyperlink}[2]{#2}
\else
\newcommand{\hrlb}{\\}
\fi

\title{The mixing time of the  giant component of a random graph}
\author{Itai Benjamini \\
{\small Weizmann Institute, Rehovot, 76100, Israel}\\
{\small {\tt itai.benjamini@weizmann.ac.il}}\\
\and Gady Kozma \\
{\small Weizmann Institute, Rehovot, 76100, Israel}\\
{\small {\tt gady\makebox[1.2 ex]{.}kozma@weizmann\makebox[1.2 ex]{.}ac\makebox[1.2 ex]{.}il}}\\
\and   Nicholas Wormald\thanks{Research supported by the
  Canada Research Chairs
Program and NSERC.}\\ {\small  University of Waterloo, Waterloo ON,
Canada N2L 3G1}\\ {\small {\tt  nwormald@uwaterloo.ca }}  }
\date{}

\begin{document}
 \maketitle
 \abstract{We show that the total variation mixing time of the simple random
 walk on the giant component
 of supercritical $\G(n,p)$ and $\G(n,m)$ is $\Theta(\log^2 n)$.  This statement was only recently proved, independently, by Fountoulakis and Reed.  Our proof  follows
 from a structure result for these graphs which is interesting in its own
 right. We show that these graphs are ``decorated expanders'' --- an expander
 glued to graphs whose size has constant expectation and exponential tail, and
 such that each vertex in the expander is glued to no more than a constant
 number of decorations.}
 \section{Introduction}
The mixing time $T$ of a finite connected graph $G$, loosely defined as the time a
random walk on that graph needs in order to be quite close to its stationary
distribution, is an important concept in randomized algorithms and theoretical
probability (see Section~\ref{mixing} for precise definitions). It has strong connections to the geometry of the graph --- for
example, $c/\lambda \leq T\leq C\log(|G|)/\lambda$ where $\lambda$ is the
spectral gap of the graph and $c$ and $C$ are some universal constants. This logarithmic factor is quite important,
though. Typically graphs coming from applications in algorithms and
statistical physics are exponentially large and have large spectral gaps, so
this factor is crucial for the applicability of an algorithm. Thus extensive
efforts have gone into understanding this factor better. See e.g.~Lov\'asz and
Kannan \cite{LK99}; Morris and Peres \cite{MP05}; or Goel, Montenegro and
Tetali \cite{GMT06}.

The mixing time of random graphs in particular is  a topic of research. Random
$d$-regular graphs were the first examples of expanders (hence they have mixing
time $\log(|G|)$). The first author and Mossel~\cite{BM}  considered the mixing time for the
simple random walk on the largest percolation cluster in a box subset of the
$d$-dimensional integer lattice. Our purpose here is to do the same for the
--- {\em a priori} one might assume simpler --- case of the largest cluster of random
graphs.

In this article, we consider two standard random graph models: the Erd\H{o}s-R\'enyi graph $\G(n,p)$ in which every edge is taken independently with
probability $p$ and $p=c/n$ for some constant $c>1$; and $\G(n,m)$, the model where all
graphs with $n$ vertices and $m$ edges are equally likely, for $m\sim cn/2$. It is well known that for such $c$ there is asymptotically almost surely (a.a.s.) a unique giant component in the random graph. Our main result is the following. Here $\Theta(f(n))$ denotes a function that, for some positive constants $c_1$ and $c_2$, lies between $c_1f(n)$ and $c_2f(n)$ for all $n$ sufficiently large.
\begin{thm}\lab{main}
Let $c>1$ and $m\sim cn/2$. The mixing time of the simple random walk on the giant component of $\G(n,m)$ is a.a.s.\ $\Theta(\log^2 n)$.
\end{thm}
It is easy to see that this implies the corresponding result for $\G(n,p)$.

The lower bound is easy, and one way to see it is this. It is straightforward to show using standard techniques that for some $c'>0$ and with $m$ as in the theorem, $\G(n,m)$ a.a.s.\ has a path of degree 2 vertices of length at least $c' \log n$. On the other hand, the hitting time of one end of a path of length $k$ is $\Theta(k^2)$.

In Section~\ref{mixing} we give a new mixing time bound (Theorem~\ref{mixthm})
for ``decorated expanders'', namely graphs which contain an expander subgraph
$B$ (we call $B$ the ``strong core'') whose deletion leaves ``small"
components such that a bounded number are
attached to any vertex of $B$. Then in Section~\ref{s:stripping} we show
(Theorem~\ref{maincore}) that $\G(n,m)$ is in fact a decorated expander. The
upper bound in Theorem~\ref{main} follows immediately from these two
results.

The construction of the strong core is quite involved. We must stress that one
cannot just take, say, the 3-core. First, it does not exist for
all $c$ down to the critical value 1. But even when it does exist, it does not
satisfy all needed properties.  See the definition of an AN-graph in Section~\ref{AN}. As part of our argument, we give an explicit result on the expansion
of random graphs with given degree sequences and all degrees at least 3 (lemma
\ref{kernel-expander}). This generalizes a number of known and folk results
--- the $d$-regular graph, the $k$-core of $\G(n,p)$ and $\G(n,m)$, the
$k$-core of bernoulli percolation on a $d$-regular graph (see \cite{G01,GM03})
and of course our strong core are all of this form.

Independently of this work, Fountoulakis and Reed~\cite{FR1,FR2}   recently
obtained the first proof of the $O(\log^2n)$ upper bound on mixing
time.   They show that the constant implicit in the bound is $O ( c^{-2})$ for $c$ bounded above by approximately $\sqrt{\ln n \ln\ln n}$.  However, their computations only relate to $\G(n,p)$, and so do not
obviously imply anything for  $\G(n,m)$ since the mixing time is not a
monotonic or even convex function on the lattice of subsets of the edge set of
the complete graph. It is interesting to compare the two
approaches. Fountoulakis and Reed start from the Lov\'asz-Kannan
integral. Applied directly it gives $C\log^3 n$. This is necessary since, as
Morris and Peres \cite{MP05} discovered the Lov\'asz-Kannan integral in fact
bounds the mixing time in the $L^\infty$ norm which in our case is really
$\Theta(\log^3 n)$. See \cite{MP05} for a detailed discussion of the
difference between these two notions. A similar problem exists with using the
spectral profile \cite{GMT06} on this
problem. Fountoulakis and Reed find a variation of the Lov\'asz-Kannan integral
which bounds only the usual, total-variation norm mixing time and gives the
correct value for the case of $\G(n,p)$. Our approach is more geometric. Once
one has that $\G(n,m)$ is a decorated expander, the mixing time is evaluated
easily with the help of the Lov\'asz-Winkler \cite{LW97} theory of
equivalences of the mixing time, a theory whose vast potential is yet to be
exploited.

Let us close this section with a few remarks

1. Consider the \emph{critical} case $c=1$. It is well
known that in this case the largest cluster is of the order of $n^{2/3}$ and
there are more clusters of comparable size. Further, it is also known that
the cluster can be split to two pieces, both of size at least $cn^{2/3}$ which are
joined by a bounded number of edges. Hence one gets that the mixing time
$T$ satisfies $T>cn^{2/3}$, or in other words, the walk does not mix rapidly at
all. This bound does not seem to be exact, though. To understand why, examine
the anomalous diffusion coefficient $\beta$. In analogy with the results for random walk on
Galton-Watson trees (\cite{BK}, see also \cite{BJKS}) one would expect that
$\beta=8/3$ or in other words that
random walk on a critical cluster would exhibit sub-gaussian diffusion. Since
the diameter of the cluster is $n^{1/3}$, the walk needs
$\left(n^{1/3}\right)^\beta=n^{8/9}$ steps to get from one end of the cluster
to the other, so a natural conjecture is that $T>cn^{8/9}$. Whether or not
this is precise, we do not know.

2. Here is an easy corollary of our structural results, for the diameter of the giant component.  The proof assumes familiarity with the later parts of the paper, but should be understandable at this point.

\begin{thm}\lab{diam}
Let $c>1$ and $m\sim cn/2$, and  $\eps>0$. The diameter of the giant component of $\G(n,m)$ is a.a.s.\ $\Theta(\log n)$.
\end{thm}

\proof
The lower bound follows easily from the property that the giant component contains induced paths of length at least $\eps \log n$ a.a.s.
The upper bound follows in two steps. Firstly,  the $\alpha$-strong core $B$, that exists a.a.s.\ in the giant component by Theorem~\ref{maincore}, is an $\alpha$-expander and consequently  has diameter $O(\log n)$.  Secondly, the  attachments to it have size, and hence height, at most $O(\log n)$ a.a.s.
\qed

We note that sharper results on the diameter of these graphs have only recently been obtained by Fernholz and Ramachandran~\cite{FR0} and Bollob{\'a}s, Janson and Riordan~\cite{BJR}.

 3. It would be interesting to extend the present result to the study of percolation on random regular graphs. It was shown in~\cite{ABS} that if the edges of a random $d$-regular graph are deleted independently with probability $1-p$ each, the threshold of appearance of a giant component is at $p=1/(d-1)$.
We conjecture that the mixing time of the random walk on the giant component for fixed $p>1/(d-1)$ is again $O(\log^2 n)$. This seems  related to the analogous question for random graphs with a given degree sequence, which would   also be of interest.

\section{Mixing times}\lab{mixing}

\subsection{Definitions}

There are many possible definition of a ``mixing time'', and extensive
literature devoted to proving relationships between the
various definitions. We shall briefly sketch the terms we shall need,
and refer the reader to~\cite{LW97} for a more orderly and far more
exhaustive introduction.
\smallskip

\noindent
{\bf Definition.\ }
Let $\mathcal{F}_0\subset\mathcal{F}_1\subset\dotsb $ be a series of
$\sigma$-fields on
a space $\Omega$. A \textbf{stopping rule} (for $\mathcal{F}_{n}$)
$\Gamma$ is a stopping time with possible external randomization,
namely, there exists some $\Omega_{2}$ such that $\Gamma:\Omega\times\Omega_{2}\to\mathbb{Z}^+$
and such that for every $\omega\in\Omega_{2}$, $\Gamma(\cdot,\omega)$
is a stopping time for $\mathcal{F}_{n}$.
\smallskip

Let $G$ be a connected finite graph, let $\sigma$ be a distribution
on the vertices of $G$, and let $\Gamma$ be a stopping rule for
a random walk $R$ on $G$ whose starting point $R(0)$ is distributed
like $\sigma$. Denote \[
\sigma_{v}^{\Gamma}=\mathbb{P}(R(\Gamma)=v).\]
In other words $\sigma^{\Gamma}$ is the distribution of the location of $R$ at
time $\Gamma$ when starting from $\sigma$. We also say that
$\Gamma$ is a {}``stopping rule from $\sigma$ to $\sigma^{\Gamma}$''.
\smallskip

 \noindent
{\bf Definition.\ }
Let $G$ be a connected graph, and let $\sigma$ and $\tau$ be two
distributions on the vertices of $G$. Define the \textbf{access time}
from $\sigma$ to $\tau$, denoted by $\mathcal{H}(\sigma,\tau)$,
using\[
\mathcal{H}(\sigma,\tau):=\min_{\Gamma:\sigma^{\Gamma}=\tau}\mathbb{E}\Gamma.\]
\smallskip

The set of stopping rules from $\sigma$ to $\tau$ is never empty:
for example, it always contains the naive rule, namely, initially
choose a vertex $v$ using $\tau$, then walk until the first time
$v$ is hit. Naturally, in most cases this rule is not optimal.

\smallskip

\noindent
{\bf Definition.\ }
Let $G$ be a connected, graph. The \textbf{mixing time} \textbf{of
$G$ } is defined as \[
\mathcal{H}:=\max_{\sigma}\mathcal{H}(\sigma,\pi)\]
where the maximum is taken over all distributions $\sigma$, and where
$\pi$ is the stationary distribution of $G$.
\smallskip

Note that we define the stationary distribution as the limit \[
\pi_{v}:=\lim_{t\to\infty}\frac{1}{t}\sum_{k=1}^{t}\mathbb{P}(R(k)=v)\]
and hence we do not need to assume that $G$ is aperiodic. Recall
that $\pi_{v}$ is proportional to the degree of $v$, $d_{v}$ namely
\begin{equation}
\pi_{v}=\frac{d_{v}}{2E(G)}\label{eq:pivdv}\end{equation}
where $E(G)$ is the number of edges of $G$.

It is interesting to compare $\mathcal{H}$ with more natural notions
of the mixing time. For example, is it true that after $t=\left\lceil C\mathcal{H}\right\rceil $
moves of a random walk, that the distribution of $R(t)$ is close
to $\pi$ in some norm? Generally the answer is no. For example, if
$G$ is a complete bipartite graph of size $2n$ with one edge added
(so that $G$ would be aperiodic). Then $\mathcal{H}\leq C$, but
it takes approximately $n^{2}$ steps until the walk becomes mixed
in the naive sense, since it needs a reasonable probability to traverse
the only edge which makes $G$ aperiodic. In this particular case,
however, it is still possible to get a uniform distribution by randomly
picking a fixed length, for example $1$ or $2$ with probability
$\frac{1}{2}$, so for practical purposes, namely for an efficient
algorithm to pick an approximately random point, it is quite reasonable
to claim that the mixing time of the graph $G$ is $\frac{3}{2}$.
It turns out that this example is typical, in the sense that by picking
the length of the walk randomly, independently of the actual steps
taken, we get close to $\pi$ by $\left\lceil C\mathcal{H}\right\rceil $
moves. For example one might take the length uniform between $1$
and $\left\lceil C\mathcal{H}\right\rceil $:
\smallskip

\noindent
{\bf Definition.\ }
Let $G$ be a connected graph. The \textbf{approximate uniform mixing
time} \textbf{of $G$} is defined by \[
\mathcal{U}_{\epsilon}:=\min_{t}||\sigma^{T}-\pi||\leq\epsilon\]
where $T$ is a stopping time with probability $1/t$ for every time
between 0 and $t-1$ independently of the walk; where $||\cdot||$
stands for the $L^{1}$ norm (a.k.a.~the total variation norm), i.e.~$||\mu-\tau||:=\sum_{v}|\mu(v)-\tau(v)|$;
and where $\pi$ is the stationary distribution of $G$ and $\epsilon>0$
is some parameter.
\smallskip

There are other variation on this {}``random number of steps'' theme.
See e.g.~\cite{A82} for a continuous time random walk version, \cite[theorem 7.2]{LW95}
for another version, and open problem 17 in chapter 4 of \cite{AF}.
\smallskip

\noindent
{\bf Definition.\ }
Let $G$ be a connected, graph. The \textbf{approximate forget time}
\textbf{of $G$} is defined by \[
\mathcal{F}_{\epsilon}:=\min_{\tau}\max_{\sigma}\min_{\mu:||\mu-\tau||\leq\epsilon}\mathcal{H}(\sigma,\mu).\]
where $\epsilon\geq0$ is some parameter, and where $||\cdot||$ stands
for the $L^{1}$ norm, i.e.~$||\mu-\tau||:=\sum_{v}|\mu(v)-\tau(v)|$.
\smallskip

$\mathcal{F}_{\epsilon}$ is called a {}``forget time'' because
we consider stopping at $\tau$ to be {}``forgetting'' our initial
distribution $\sigma$. Perhaps surprisingly, the minimum $\tau$
is not necessarily achieved at $\pi$, and for directed graphs the
ratio $\mathcal{H}/\mathcal{F}_{0}$ may be arbitrarily large. See
\cite{LW97} for a detailed discussion.

\begin{thm}
\label{thm:lovasz}For any $\epsilon\leq\frac{1}{2}$, $\mathcal{F}_{\epsilon}=\Theta(\mathcal{H})=\Theta(\mathcal{U_{\epsilon}})$
where the constants implicit in both $\Theta$ may depend on $\epsilon$.
\end{thm}
The inequality $\mathcal{F}_{\epsilon}=\Theta(\mathcal{H})$ follows
from theorems 3.1, 3.2 and 3.8 in \cite{LW97} --- note that a random
walk on a (non-directed) graph is always a time-reversible Markov
chain. The inequality $\mathcal{U}_{\epsilon}=\Theta(\mathcal{H})$
comes from corollary 5.4 ibid.
\smallskip

\noindent
{\bf Definition.\ }
The (edgewise) Cheeger constant of a connected graph $G$ is defined
by \[
\Phi:=\min_{0<\pi(S)\leq\frac{1}{2}}\frac{1}{\pi(S)}\sum_{i\in S,j\not\in S}\pi_{i}p_{ij}\]
where $p_{ij}$ is the probability to step from $i$ to $j$, namely
$1/d_{i}$ if $j$ is a neighbour of $i$ and $0$ otherwise; where
$\pi$ is the stationary distribution of $G$, and where $\pi(S):=\sum_{i\in S}\pi_{i}$.
\smallskip

Note that plugging in (\ref{eq:pivdv}) we see that the element inside
the $\min$ is, more or less the quotient of the number of edges leading
out of $S$ divided by the number of edges inside $S$, hence the name
{}``edgewise'' Cheeger constant. A similar value is called {}``conductance''
in \cite{JS89,LK99}.

We shall need the following connection between the Cheeger constant
and the mixing time.

\begin{thm}
\label{thm:JS}Let $G$ be a connected aperiodic graph. Then\[
\mathcal{H}\leq C\log(1/\min\pi_{i})\frac{1}{\Phi^{2}}.\]

\end{thm}
This was first proved by Jerrum and Sinclair in \cite{JS89} (the particular
case of expanders, which is what we will use, was proved earlier by
Alon \cite{A86}, and in the continuous setting this goes back to Cheeger
\cite{C70}). We shall
only use the form
\begin{equation}
\mathcal{H}\leq C\log   E(G) \frac{1}{\Phi^{2}}\label{eq:JSE}\end{equation}
which follows immediately from (\ref{eq:pivdv}).

\subsection{AN-graphs.}\lab{AN}

\noindent
{\bf Definition.\ }
We say that a connected graph $G$ is an $\alpha$\textbf{-AN graph},
(or an $\alpha$-decorated expander) where $\alpha>0$ is some number,
if the graph has a subgraph $B$ with the following properties:
\begin{enumerate}
\item $B$ is a $\alpha$-expander, i.e.~$\Phi(B)\geq\alpha$.
\item The connected components $D_{i}$ of $G\setminus B$ are small in
the following sense: denote by $E'(D_{i})$ the number of edges in
$G$ with at least one vertex in $D_{i}$, or in other words, the
internal edges of $D_{i}$ added to the edges connecting $D_{i}$
to $B$. Then\begin{equation}
\#\{ i:E'(D_{i})\geq\lambda\}\leq   E(G) e^{-\lambda\alpha}.\label{eq:expotails}\end{equation}
In particular there are no components with $E'(D_{i})>\frac{1}{\alpha}\log E(G)$.
\item Each $v\in B$ is connected to no more than $\frac{1}{\alpha}$ different
$D_{i}$-s.
\end{enumerate}
\smallskip

$G\setminus B$ here denotes the graph reached after removing the
vertices of $B$ and all edges with at least one vertex in $B$ from
$G$. Note that the definition is meaningless for $\alpha>1$, so
we will always assume $\alpha\leq1$.

\begin{thm}\lab{mixthm}
The mixing time of an $\alpha$-AN graph is $\leq C\alpha^{-6}\log^{2} E(G) $.
\end{thm}
Note that the Cheeger constant of $G$ might be $\leq C\alpha^{3}/\log E(G) $.
For example, take a subset of $B$ realizing the maximum --- assume
it is small --- and hang from each vertex $\frac{1}{2\alpha}$ copies
of straight line segments of length $\log  E(G)/2\alpha$. Hence theorem
\ref{thm:JS} only gives a bound of $C\alpha^{-6}\log^{3} E(G)$. In
our application $\alpha$ will be a constant independent of $n$,
so the result of the theorem is an asymptotic improvement. This example
also shows that the mixing time is $>c\alpha^{-2}\log^{2} E(G)$ since
if we start from the end of one such straight line segment we need
that many steps to have a decent probability to exit it.

\begin{proof}
For every $v,w\in B$, let $q(v,w)$ be the probability that a random
walk on $G$ starting from $v$ hits $B$ in $w$ (in particular,
$q(v,w)\neq0$ only if $v$ and $w$ are   neighbours in $G$ or if $v$
and $w$ have neighbours in the same $D_{i}$). The symmetry of the
random walk on $G$ gives $d_{v}(G)q(v,w)=d_{w}(G)q(w,v)$ --- we
will denote the degree of a vertex by $d_{v}(G)$ when it is not clear
about which graph we are talking. Now construct a weighted graph with
self loops $B_{1}$ with the vertex set identical to the vertex set
of $B$, and for any $v,w\in B_{1}$ make the weight of the edge between
$v$ and $w$ be $d_{v}(G)q(v,w)$. It is easy to see that $B_{1}$
is $B$ with some added edges and some edges with increased weight,
and the total weight of every vertex is increased by at most $\frac{1}{\alpha}+1\le \frac{2}{\alpha}$
due to requirement 3 from an $\alpha$-AN graph. Hence \begin{equation}
\Phi(B_{1})\geq\tfrac 12\alpha^{2}.\label{eq:PhiH1}\end{equation}
We now use the clause $\mathcal{H}=\Theta(\mathcal{U}_{\epsilon})$
in Theorem \ref{thm:lovasz}: define $\pi_{1}$ to be the stationary
distribution of $B_{1}$; let $C_{2}$ be some constant sufficiently
large and define \begin{equation}
J=\left\lceil C_{2}\mathcal{H}(B_{1})\right\rceil .\label{eq:defJ}\end{equation}
Let $S_{1}$ be the stopping rule stopping at time $t=0,
\ldots
,J-1$
with probability $\frac{1}{J}$. Let $\sigma_{1}$ be any distribution
on $B_{1}$. Then we get, for $C_{2}$ sufficiently large,\begin{equation}
||\sigma_{1}^{S_{1}}-\pi_{1}||\leq\frac{1}{4}\quad\forall\sigma_{1}.\label{eq:sigma1S1}\end{equation}
Also, Theorem \ref{thm:JS} allows us to estimate \begin{equation}
J\leq\left\lceil C\log E(B_{1})\Phi^{-2}(H_{1})\right\rceil \stackrel{(\ref{eq:PhiH1})}{\leq}C\alpha^{-4}\log E(B_{1})\leq C\alpha^{-4}\log E(G).\label{eq:estJ}\end{equation}

Next define stopping times $\tau_{j}$ as follows: $\tau_{0}:=-1$
and\[
\tau_{j}:=\min\{ t>\tau_{j-1}:R(t)\in B\}.\]
The point about the definition of $B_{1}$ is that the regular random
walk on $B_{1}$ is identical to the process $R(\tau_{j})$. This
allows to translate (\ref{eq:sigma1S1}) to the setting of the random
walk on $G$: define $S$ to be the stopping time that stops at $\tau_{j}$,
$j=1,
\ldots
,J$ with probability $\frac{1}{J}$. Then for any starting
distribution $\sigma$ on $G$,\[
||\sigma^{S}-\pi_{1}||\leq\frac{1}{4}.\]
Reaching the distribution $\pi_{1}$ is in effect forgetting $\sigma$,
hence we may use the clause $\mathcal{H}=\Theta(\mathcal{F}_{\epsilon})$
of Theorem~\ref{thm:lovasz} to get that \[
\mathcal{H}\leq C\max_{\sigma}\mathbb{E}S\leq C\max_{v}\mathbb{E}\tau_{J}.\]

Hence we have the task of estimating $\mathbb{E}\tau_{j}$. In general,
if $D$ is any graph and $v\in D$; and if $R$ is a random walk starting
from $v$, then the expected time until $R$ returns to $v$ is $1/\pi(D)_{v}$.
See \cite[chapter 2, lemma 5]{AF}. Take some $j>0$ and let $D\subset G$
be 
$$
D:   =\bigcup\{ D_{i}:D_{i}\textrm{ neighbours }R(\tau_{j})\}\cup
    \{ v\in B:\exists i\textrm{ s.t.~}D_{i}\textrm{ neighbours }R(\tau_{j})\textrm{ and }v\}.
   $$
The definition of an $\alpha$-AN graph shows that there are $\leq\alpha^{-1}$
$D_{i}$-s neighbouring $R(\tau_{j})$ and each one satisfies $E'(D_{i})\leq\alpha^{-1}\log E(G)$,
so $E(D)\leq\alpha^{-2}\log E(G)$, and then\begin{equation}
\mathbb{E}(\tau_{j+1}-\tau_{j})\leq1/\pi(D)_{R(\tau_{j})}\stackrel{(\ref{eq:pivdv})}{\leq}2E(D)\leq2\alpha^{-2}\log E(G).\label{eq:taujp1tauj}\end{equation}
This estimate does not work for $\tau_{1}$ since $R(0)$ might not
belong to $B$. Here we need the fact that for any graph $D$, the
expected time that $R$ takes from $v$ to $w$ is $\leq2E(D)\cdot\rho$
where $\rho$ is the electrical resistance between $v$ and $w$.
See \cite[chapter 3, corollary 11]{AF}. Clearly $\rho\leq|D|$ since
the resistance between $v$ and $w$ is at most the resistance of a
path betwen them. Let therefore $D_{i}$ be the component containing
$R(0)$ and let $D$ be $D_{i}\cup\{\textrm{its neighbours in }B\}$.
We get \begin{equation}
\mathbb{E}\tau_{1}\leq E(D)^{2}=E'(D_{i})^{2}\leq\alpha^{-2}\log^{2}E(G).\label{eq:tau1}\end{equation}

Collecting (\ref{eq:estJ}), (\ref{eq:taujp1tauj}), and (\ref{eq:tau1})
we get that \[
\mathbb{E}\tau_{J}\leq C\alpha^{-6}\log^{2}E(G)\]
which finishes the proof. \qed
\end{proof}

 \section{Random graph preliminaries}~\lab{prelim}
We define an  {\em  $\alpha$-strong core} of a graph $G$ to be any  subgraph $B$ of with the properties as listed in the definition of the $\alpha$-AN graph in Section~\ref{AN}.
 We seek an $\alpha$-strong core of
$G\in\G(n,p)$ for $p=c/n$ where
$c>1$ is fixed. It suffices (and indeed gives a stronger result) to consider
 $G\in\G(n,m)$ for $m\sim cn/2$. (See Bollob{\'a}s~\cite{BB} or Janson et
 al.~\cite{JLR} for these basic definitions and results on random graphs).
 The 2-core of a graph  or multigraph is the maximum subgraph of minimum
degree at least 2. The 2-core, if it exists, is known to be unique;
otherwise we say the 2-core is empty. It can be obtained by recursively
deleting vertices of degree 0 and 1. Define
$$
b=b(c):=1-t/c,
$$
where $t=t(c)$ is the unique root of the equation
$$
te^{-t}=ce^{-c},\quad t\in (0,1).
$$
The following results are well known. (See~\cite{JLR} for example,
or~\cite{PW} for more precise results on the joint distributions and related
information. See~\cite{Whandbook} for definition of the combination of the
notations $o$ and a.a.s.)
\begin{thm}\lab{giantconc}
  The number of vertices of the giant component of
$G\in\G(n,p)$ is  a.a.s.\ $bn+o(n)$, the number of vertices of the
2-core of $G$   is a.a.s.\ $b(1-t)n+o(n)$, and  the
number of edges of the 2-core is a.a.s.\ $b_1n+o(n)$ where $b_1=b_1(c)>b(1-t)$. The same results hold for $G\in\G(n,m)$ with $m\sim cn/2$.
\end{thm}
 Deleting the edges of the 2-core   therefore a.a.s.\ leaves a forest $F$ of
$b(1-t)n+o(n)$ trees  with   $bn+o(n)$   vertices in total.
  We will condition on the event that the sizes of the giant and 2-core satisfy
these conditions. Each forest of a given number of trees with a given
number of  vertices is equally likely to occur as $F$. Also well known is the following type of result. First,  we say that the distribution of a random variable $X$ has an {\em exponential tail} if $\pr(X>j)=O(e^{-cj})$ for some $c>0$.
\begin{lemma}\lab{treesize}
Let $g(n)$ be a fixed function with $g(n)=o(n)$. Let  $G\in\G(n,p)$, conditional upon the giant component having between $bn-g(n)$ and $bn+g(n)$ vertices, and the 2-core having between $b(1-t)n-g(n)$ and $b(1-t)n-g(n)$ vertices. The size of the tree in $F$ containing a given vertex of the 2-core (conditional upon that vertex being in the 2-core) has an exponential tail.
\end{lemma}
\proof
We may fix the 2-core with $s\sim b(1-t)n$ vertices, and assume the vertices
of the giant componenet not in  the 2-core are labelled $1,\ldots,r$ (where
$r\sim btn$). Then, after deleting all edges of the 2-core, each forest of $s$   trees with  root vertices in the 2-core (mutually distinguishable from each other but unlabelled), and the non-root vertices labelled $1,\ldots,r$, is equally likely to occur  as $F$. The number of such forests is $s(r+s)^{r-1}$ (see~\cite[p. 17]{M}). The number of possible trees rooted at the first 2-core vertex, given the tree has $j+1$ vertices (that is, $j$ labelled vertices plus the root vertex) is ${r \choose j}(j+1)^{(j-1)}$, where the first factor chooses the tree's vertices and the second constructs the tree. Additionally, of course, the remaining part of the forest is counted by $(s-1)(r-j+s-1)^{r-j-1}$. Applying Stirling's formula and a little manipulation now shows that the probability that the first tree has size $j$ is
$$
O\big(sj^{-3/2}\big((1-\rho)e^\rho)\big)^j\big)
$$
where $\rho=  s/(s+r)$. Since $c$ is fixed, $(1-\rho)e^\rho$ is less than, and bounded away from, 1. The lemma follows.
\qed

The 2-core will not usually be an $\alpha$-strong core  of the giant component
because it has  long (length $c' \log n$) paths of degree 2 vertices.  The
$\alpha$-strong core we will be using is obtained from the 2-core by deleting paths
of degree 2 vertices. We make this precise as follows. Define a {\em
2-path} in a graph $G$ to be a path induced by vertices of degree 2 in $G$, and   an {\em
isolated cycle}  to be a component of $G$ that is just a cycle.

Analogous to the definition for a random variable, we say that a  set  $S=\{s_i\}$ of nonnegative numbers  has an
exponential tail if some $C>0$ and $C'$ exist such that
 for all
$j\ge 0$
$$
\frac{|\{ i :  s_i\ge j\}|}{|S|} <C' e^{-Cj}.
$$
(Equivalently, we could work with the definition in which $C'=2$.)
Moreover, if the  set $S$ is indexed by $n$ and there exist universal constants $C$ and $C'$ for
which the inequality is true a.a.s. (as $n\to\infty$), we say that $S$ has
an {\em exponential tail a.a.s. }

We obtain the following result quite easily using Markov's inequality and  Lemma~\ref{treesize}, together with simple sharp concentration  of the numbers of trees of a fixed size (say by applying Chebyshev's inequality after computing second moments along the lines of the calculation in Lemma~\ref{treesize}).  The proof is left as an exercise.
\begin{lemma}\lab{treesizes}
The set of sizes of trees in $F$ a.a.s.\ has an exponential tail.
\end{lemma}

\section{Stripping processes} \lab{s:stripping}
 Let $N$ be any integer. Given a graph, we can perform a ``stripping"
process that first removes all vertices  not in the 2-core, $G$, of the graph and then  recursively
deletes isolated cycles and the vertices in 2-paths of length greater than
$N$, as well as vertices of degree less than 2 (which arise if both ends of
a 2-path are adjacent to a vertex of degree 3)   and vertices with more
than $N$ removed neighbours, until no more deletions are
possible. Note that the resulting graph does not depend on the choices made
at each step, since, once a vertex can be deleted, it remains deletable in
any subsequent step. Presumably when $G$ is   the  2-core of the giant
component of a graph in $ \G(n,m)$ (with $m$ as in Theorem~\ref{giantconc}), this resulting graph is, for $N$ sufficiently large,   an $\alpha$-strong  core a.a.s. However, to make the proof easier we will modify
this stripping process.

To analyse such  processes we will consider the {\em kernel} $K(G)$ of
$G$, defined for $G$  with minimum degree at least 2. The kernel is
obtained by replacing each maximal 2-path that joins two vertices $u$ and
$v$ of degree at least 3 by an edge $uv$,  and deleting each isolated
cycle. It is possible that the kernel possesses loops and/or multiple edges.
A loop contributes 2 to the degree of its incident vertex, so
$\delta(K(G))\ge 3$ ($\delta$ denoting minimum degree).

The kernel was used in~\cite{PW} to derive properties of the 2-core of the
random
$G\in\G(n,m)$. In particular, it is easy  to obtain the following from the results there.
Here $b$ is the same as defined in Section~\ref{prelim}.
\begin{lemma}\lab{kernelconc}
The number  of
vertices of degree 2 in the 2-core of $G\in\G(n,m)$ is  a.a.s.\
$ b_2n+o(n)$ for a constant
$b_2$,   depending on $c$, with $0<b_2<b(1-t)$.
\end{lemma}
 It follows from this and Theorem~\ref{giantconc} that the
size (number of vertices)
 of the kernel is similarly sharply concentrated at $(b-b_2)n$, and the   number  of edges in the kernel is also sharply concentrated at $(b_1-b_2)n$.

In examining the stripping process, it is difficult to keep track of the
distribution of lengths of those 2-paths containing vertices of degree 2
that still remain but were adjacent to removed vertices. So we  define
another stripping process, called {\em severe stripping}, that in general
removes more  than is necessary, always erring on the safe side. This can be
applied to a graph  or multigraph $G_0$ with $\delta(G_0)\ge 2$ and with no isolated cycles. To guide the process,   some edges and vertices of the kernel are designated as red. Colouring an edge or vertex of the kernel of a graph red marks the corresponding part of the graph for removal during the severe stripping process. All edges incident with red vertices are also painted red.

To initialise this process, begin with any graph $G$ and obtain $G_0$ from the 2-core of $G$ by deleting the isolated cycles.  For any vertex $v$ of $G_0$ that is adjacent to at least $N-1$ vertices of $G-G_0$, if $d_{G_0}(v)=2$ then the edge of the kernel $K(G_0)$ corresponding to the maximal 2-path containing $v$ is painted red, whilst if  $d_{G_0}(v)\ge 3$ then it is a vertex of
$K(G_0)$ and is painted red, as are all the edges of $K(G_0)$ incident to $v$.  Also, with $N$ as above, take
note of all maximal 2-paths of $G_0$ containing more than $N/2$ vertices:  the  corresponding edges of $K(G_0)$ are  coloured red too. For later reference, we call $G_0$ the {\em trimmed core} of $G$, and the graph $K(G_0)$ together with its colouring defined in this way is called the {\em painted kernel} of $G_0$ with respect to $G$.

 During the process, some other kernel edges will be coloured purple, and of course the kernel will be modified as the graph changes. Purple edges correspond  to 2-paths whose lengths have been ``exposed"  in the sense
that they have been used to influence the algorithm, and yet which are not (yet) required to be removed. Also, some kernel vertices will be coloured pink, to signify that they have already lost a neighbour during the process. Any pink vertex that loses a second neighbour is immediately recoloured red.  This is  to ensure, without excessive bookkeeping, that, by the end of the process, all remaining vertices have lost at most $N$ neighbours each. (We may assume $N\ge 2$.) The looseness this causes in the final bound only affects terms that we are not attempting to maximise.

To simplify the argument further, we
will also avoid keeping a record of the length of a purple edge at any steps
after it is first formed. Hence, we must remove a purple edge whenever
it is merged to another edge, i.e.~when a vertex $v$ adjacent to one of its
ends drops to degree 2.

Formally, severe stripping defines a sequence $G_0, G_1, \ldots$ of
graphs and colourings of their kernels as
follows, beginning with the trimmed core $G_0$ and its painted kernel with respect to $G$.

For step $i$, select  a red edge   of
$K(G_{i-1})$ uniformly at random from all the red edges. Remove the corresponding maximal
2-path $\gamma_0$ of $G_{i-1}$, to obtain the next graph $G_i$.   Let $u$ be the
vertex of $G_i$ adjacent in $G_{i-1}$ to one end of $\gamma_0$.  If $u$ still has degree at least 3, colour it pink if it was uncoloured in $K(G_{i-1})$, whilst if $u$ was already pink,  colour $u$ and all  edges of $K(G_{i})$ incident with $u$ red.  On the other hand, if $u$ has degree 2 in $G_i$, let $e$ be the   edge of
$K(G_i)$ corresponding to the maximal 2-path $\gamma_1$ of $G_i$ containing $u$.
Colour $e$ red if
either of the two edges of $K(G_{i-1})$ which form $u$ was already
purple  or red, and colour $e$ purple otherwise. All other edges and vertices of
$K(G_i)$ inherit their colours from
$K(G_{i-1})$.   All red vertices will be deleted eventually --- when
their degree drops to two they will be deleted from the kernel and will remain
as part of the path of the graph $G_i$ corresponding to a red edge of the
kernel, and when that edge is deleted the vertex will be removed from $G_i$ (in
some special cases below, both deletion steps will happen at the same time).

 There is a special case: if $u$   is contained in an isolated cycle of $G_i$, then the whole cycle is removed.  Treat the vertex $u'$
adjacent to the other end of
$\gamma_0$  separately with the same rules. If $u=u'$  the rules above apply in
the obvious way, unless the degree of
$u$   falls to 1 when $\gamma_0$ is removed. In this  case, remove the maximal
path of degree 1 and 2 vertices that contains $u$ from $G_{i-1}$, let $u''$
be its point of attachment and repeat the above colouring rules treating
$u''$ as $u$.

The severe stripping process continues, repeating the above step, until the
point is reached that no red edges remain in $K(G_M)$.  The  {\em
$N$-reduced core}  of a graph $G$, denoted $R_N(G)$, is the final graph $G_M$ obtained by applying the process starting with the trimmed core and its painted kernel with respect to $G$. This may seem to depend on the order of choosing the red edges for removal, but it is actually unique, which
is convenient for descriptive purposes but unimportant for our arguments. The
uniqueness can be seen, by observing that severe stripping is equivalent to recursively removing all red edges and vertices, any vertex which has had at least two incident edges removed, and any edge of the kernel that at any point has more than $N$ degree 2 vertices or comes from merging at least three of the edges of $K(G_0)$.

One aspect of the definition of severe stripping may seem redundant at this point: one could avoid painting vertices red if one instead painted all incident edges red. The resulting process would be the same, but in the analysis we need to know which vertices of the random graph have been investigated in some sense, so for this reason the red vertices are recorded.

Using the severe stripping process we will obtain the following.

\begin{thm}\lab{maincore}
Fix $c>1$ and let $G\in \G(n,m)$ where $m\sim cn/2$. For
$N$  sufficiently large (depending on $c$) and $\alpha$ sufficiently small, $R_N(G)$ is a.a.s.\
an $\alpha$-strong core of $G$.
\end{thm}
The proof of this theorem is spread out over the next two sections.

 For simplicity we redefine $n$ and $m$ so that $G\in\G(\hat n,\hat m)$, and we condition on the numbers $n$  of vertices and $m$ of edges in the 2-core of
$G$.  The 2-core is distributed u.a.r.\ (uniformly at random) as a graph with these parameters $n$ and $m$,  and minimum degree at least 2. The parameters are sharply
concentrated as discussed above, so we examine for a while the random graph
space  $\G_2(n,m)$ containing all $(n,m)$ graphs with minimum degree at
least 2.

To study this, and in particular the severe stripping process, we
take the approach in~\cite{CW}, which reveals that the following model, used by Bollob{\'a}s and Frieze~\cite{BF} and  Chv\'atal~\cite{C}, is very convenient
for such purposes. (A similar idea was used in~\cite{AFP}.)  A random element of this model, which we call
$\C(n,m)$, is obtained as follows. Start with
$n$ isolated vertices and add
$m$ edges by choosing each end of each edge uniformly at random. All
choices are made independently with replacement. This is equivalent to choosing $m$ labelled oriented edges and then forgetting their labels and orientations. The result is a
pseudograph that can have loops and multiple edges. Note that the
restriction to simple graphs gives precisely $\G(n,m)$ (the uniform space).

Now define the probability space $\C_k(n,m)$  to be the
  restriction of $\C(n,m)$ to the graphs with minimum degree at least $k$.
We will proceed to analyse severe stripping applied to $\C_2(n,m)$, for
appropriate values of $m$.

{From}~\cite{CW} for example, we know the following.
\begin{lemma}\lab{tocore}
If we repeatedly delete vertices of degree 0 and 1 from
$\C(n,m)$, the final result, conditional on its numbers $n'$ and $m'$  of
vertices and edges, is distributed precisely as  $\C_2(n',m')$.
\end{lemma}
The proof is   simple enough to be omitted, using induction on the
steps of the deletion process (see below for more complicated applications of
this technique). With a similar
step-by-step approach, we easily obtain the following, which bears some resemblance to~\cite[Lemma~3]{PW}.    By {\em
suppressing} a vertex
$v$ of degree 2, we mean joining its two neighbours  with a new edge and
then deleting $v$.

As {\L}uczak~\cite{L} observed, we know that  a.a.s.\  there will only be a small number of vertices in isolated cycles in   $\C_2(n,m)$, which will affect the parameters after those cycles are discarded. To avoid switching notation after deleting a small number of vertices, we proceed initially as if no isolated cycles occurred. Let $\C_2(n,m)^*$ be the probability space derived from $\C_2(n,m)$ by restricting to those pseudographs with no isolated cycles. Note that for the following lemma and similar statements, assuming that $n$-vertex graphs have vertex set $[n] =\{1,\ldots, n\}$, we should map the $n'$ vertices of the kernel into the set $[n']$ in a canonical way. Each time a vertex of degree 2 is deleted, the remaining vertices may be renumbered, preserving the ordering. This renumbering is sometimes done implicitly in our arguments.

\begin{lemma}\lab{C23}
If we begin with a random member $M$ of $\C_2(n,m)^*$ and
suppress vertices of degree 2  repeatedly
 until none remain, the result, conditional on its numbers $n'$ and $m'$  of
vertices and edges, is distributed precisely as  $\C_3(n',m')$.
\end{lemma}
\proof
For this proof, we may retain the labels and orientations of the edges of $M$ as in the definition of the model. Then, conditional on the set of degree 2 vertices in $M$, it is uniformly distributed. To suppress a degree 2 vertex, choose the vertex $v$ to be suppressed and the end of an incident edge $x$ to delete.  The other edge incident with $v$ is extended to meet the vertex at the other end of $x$, while both $x$ and $v$ are deleted, to obtain a pseudograph $M'$ (with oriented, labelled edges). It is clear that the number of ways to reverse this operation is independent of $M'$, given the labels of $v$ and $x$ (which must be missing in $M'$). By induction, the pseudograph obtained after suppressing $k$ of the degree 2 vertices in this manner is uniformly distributed, given its set of vertices, edges, and degree 2 vertices. The lemma follows from this statement applied to  $k$ being the number of vertices of degree 2 in $M$.\qed

{From} Lemma~\ref{C23}, the kernel  $K(G_0)$ can  be modelled by  $\C_3(n',m')$, where
$(n',m')$ will be restricted to the range of the sharp concentration shown
above  from known results about
  simple graphs. Our conclusions that are a.a.s.\ true for kernels with this
range  of values will then be shown to apply to the case that the initial
graph was simple.

  We call it
a  {\em random ordered assignment} of a given set of vertices to a given
set of edges of a pseudograph if the given vertices  are randomly
assigned to those edges and the ones assigned to a particular edge placed along it in some order, such that, with parallel edges   canonically
distinguished from each other and loops   given a canonical direction, each
assignment (including the ordering along each edge) is equally likely.
\begin{lemma}\lab{C32}
The  elements of $\C_2(n,m)^*$ having a kernel with vertex set $[n']$ and with $m'$ edges are obtained with uniform distribution
by starting with the kernel  randomly taken from $\C_3(n',m')$, and then
using a random ordered assignment $g$ of the vertices of degree 2 to the edges
of the kernel.
\end{lemma}
 The proof is omitted as it is very similar to the proof of Lemma~\ref{C23}. The version for simple graphs  was used by {\L}uczak~\cite{L} and in~\cite[Section~4]{PW} where
kernels are examined conditioned on degree sequence.

Lemma~\ref{C32} shows that to analyse the severe stripping algorithm applied
to  $\C_2(n,m)^*$, we may consider a random kernel taken from $\C_3(n',m')$
and  a random ordered assignment $g$ of a given set $S_0$ of degree 2 vertices
(for all the appropriate values of the parameters).  We need to use a version
of the method of deferred decisions: we do not examine the end of any purple
or red edge until it is needed for a decision in the severe stripping
algorithm.  For a   precise description, we argue as in~\cite{CW} but with a
model similar to the kernel configuration model of~\cite{PW}. Model
$\C_3(n',m')$ as the set of random functions $f$ from $[2m']$ to $[n']$. The pair of vertices $f(2j-1) f(2j)$ forms an edge for $1\le j\le m'$. Thus, the function $g$ maps the set of degree 2 vertices (of which there are $m-m'$) to $\{j:1\le j\le m'\}$, where the number $j$ represents the edge $f(2j-1) f(2j)$. We say that this edge has {\em label} $j$.

The presence of a painted kernel affects the process, so we have to define some sets for special attention.
 At the end of the $i$th step of the stripping process, let  $S_i\subseteq S$ denote the set of degree 2 vertices on noncoloured edges (i.e.\ non-red non-purple edges) of
$G_i$. Let $R_i$ denote the set of all $j$ such that $f(2j-1) f(2j)$ is a red edge of  $K(G_i)$, and $P_i$ the corresponding set for
purple edges.   Also $VR_i$ and $VP_i$ are the sets of red and pink vertices respectively.
As with the vertex labels, at each step that an edge is deleted, the edge
labels are compressed into the range $[m']$ (where $m'$ is the number of
edges of the kernel of the new graph) and the action of $f$ and $g$ is
modified accordingly. When two edges of the kernel coalesce into one (due to
a common adjacent vertex being reduced to degree 2), a similar canonical
relabelling of the edges is carried out in which the new edge is given some canonical label, say the largest edge label.

In the end we will show that we only need to deal with a stripping process with the starting graph $G_0$ drawn uniformly at random from $\C_2(n,m)^*$. We can generate    a random element of this model while performing the stripping algorithm,  ``exposing" only those parts of the graph as required for steps of the algorithm. After step $i$ the exposed parts are the labels of the purple and red edges of $K(G_i)$, the degree 2 vertices on each of these (i.e.\ that part of the ordered assignment $g$), and all the values $f(i)$ contained in $VR_i$.
Thus, initially the labels of the red vertices and their preimage under $f$, and the members of $S_0$ assigned to red edges, are all given.
When the vertex $u$ at the end of a red edge with label $j$ is investigated, the value $f(2j-1)$ or $f(2j)$, as the case may be, is first chosen from the non-red (i.e.\  uncoloured and pink) vertices, and then it is decided (randomly, with the correct probability, which the following lemma gives a simple way to calculate) if $u$ has degree 3.    If $u$ has degree greater than 3, it is coloured  pink but the remaining part of $f^{-1}(u)$ is kept random, i.e.\ not exposed. If degree 3, the two adjacent edges are determined (i.e.\ the two remaining elements of $f^{-1}(u)$ are  decided), as is the part of $g$ assigning vertices to these edges. Provided that these two edges are distinct, they coalesce into one new purple or red edge (depending on the sizes of those preimages). For this new edge, only its label $m'$  is known in this step, and not the endvertices of the edge.  On the other hand, if the two elements of $f^{-1}(u)$ belong to the same edge, it is simply removed,
 because it corresponds to the appearance of an isolated cycle in the stripping algorithm.

The following lemma asserts that the unexposed part of the graph remains
nicely random, in order for the whole ``exposing" process to work as
described. First, note that we defined the trimmed core and the
painted kernel for multigraphs and hence they apply to members of $\C(n,m)$.
\begin{lemma} \lab{strip}
Let $G_0$ be the trimmed core of a random   multigraph $G$ in $\C(n,m)$,
and colour $K(G_0)$ as   the painted kernel of $G_0$ with respect to $G$. Then
apply the severe stripping process to obtain $G_1,G_2,\ldots$. Next,
condition  on $|V(K(G_i))|= n'$  and  $|E(K(G_i))|= m'$,   on
the sets $R_i$,
$P_i$, $VR_i$ and $VP_i$, and on $f^{-1}(j)$ for all $j\in VR_i$.   Then the remaining values of $f$ are distributed uniformly at random on $[n]\setminus  VR_i$  conditional upon $|f^{-1}(j)|\ge 3$ for all $j$.
\end{lemma}
\proof
A key thing to realise for this proof is that the conditioning described does not put any constraints on the relative positions of the coloured edges and vertices.

The lemma is proved by induction on $i$.  This is   similar to the proof of Lemma~5 in~\cite{CW}, only
more complex because various cases of encountering red or purple edges need
to be considered. The case $i=0$ follows from Lemmas~\ref{tocore} and~\ref{C32} since the number of ways to reinstate the isolated cycles is
independent of the values of $f$.

We will now show how step $i+1$ follows from step $i$. Let $R_{i+1}$,
$P_{i+1}$, $VR_{i+1}$ and $VP_{i+1}$ be
given, and let $f_1$ and $f_2$ be two functions from $[2m']$ to
$[n']$ such that $|f_k^{-1}(j)|\geq 3$ that satisfy the compatibility
condition that $f_k(j)\in VR_{i+1}$ implies that $f_{3-k}(j)=f_k(j)$. We
need to show that both functions $f_k$ have the same probability. The function $f_k$ could come from
some $g_k$ in the $i$th step by various means. Let us take as an example the
case where a red edge connected to two uncoloured vertices of degree at least $4$
is removed (and the vertices are hence coloured pink). To reverse this
process, one must find two
pink vertices $p_1$ and $p_2$, uncolour them, add an edge between them at
some position $s\in[m']$, and relabel the edges. Examine some specific $p_1$,
$p_2$ and $s$. It is easy to see that they dictate the sets $R_i$, $P_i$,
$VR_i$ and $VP_i$ and further that the two functions $g_1$, $g_2$ satisfy
  that $|g_k^{-1}(j)|\geq 3$ as well as the same
compatibility condition the functions $f_k$ satisfy. Hence (by induction) they have
the same probability. Further, since the red edge to be removed is selected
randomly, the probability of $s$ to be selected for removal for $g_1$ is the
same as for $g_2$. Hence they contribute the same amount to the functions  $f_k$. Since
this holds for any values of $p_1$, $p_2$ and $s$, we get that the total
contribution of our example case (a red edge connected to two uncoloured
vertices of degree $\geq 4$) is the same to $f_1$ and to $f_2$.

As a second example, let us take the case that one end is uncoloured with
degree 3, and its other two edges are also uncoloured. The stripping rules
require us to suppress the vertex and colour the new edge purple. Hence the
reversal process consists of finding a pink vertex $p$, verifying that the
 highest labelled edge $q$ is purple, finding
a location for the removed vertex $t$, three locations for new edges, $s_1$ for
the red edge removed and $s_2$, $s_3$ for the uncoloured edges merged. As
above we see that the reversal process does not depend on the ``remaining
values of $f$'' and the argument goes through unchanged. We will not bore the
reader with any more cases.
  \qed

\bigskip

 This lemma is used  in the proof of the next result.  We will use
$\eps_N$ to denote some function that is  constant for fixed $N$
but  goes to 0   as
$N\to\infty$ (perhaps different functions at different occurrences of the
notation).
\begin{lemma} \lab{epsilons}
Let $m\sim cn/2$ with $c>1$. Let $G_0$ be the trimmed core of $G\in\C(n,m)$, and let $\hat n$ and $\hat m$ be the (random) numbers of  vertices and edges of $K(G_0)$.
  The numbers $n'$ and $m'$ of vertices and edges of  the kernel $K(R_N(G))$ of the $N$-reduced core of $G$ a.a.s.\ satisfy $n'=\hat n-O(\eps_N n)$ and $m'=\hat m-O(\eps_N n)$. Moreover, conditional upon having particular values of $n'$ and $m'$, $K(R_N(G))$  is distributed as $\C_3(n',m')$.
\end{lemma}
\proof
 The simplest part is the last, as it follows directly from lemma
\ref{strip}. Indeed, for any possible value $H$ of $K(R_N(G))$, its probability
comes from a sum over all its realizations as a function $f:[2m']\to[n']$, all
possibilities for the number of trimming steps $i$ and all possibilities for
$P_i$ and $VP_i$. However, the number of possibilities does not depend on the
structure of $H$ at all, and the probability of each quadruple $f, i, P_i,
VP_i$ does not depend on $f$ by lemma \ref{strip}. This shows that
$K(R_N(G))$, conditioned on $m'$ and $n'$ is indeed distributed as
$\C_3(n',m')$.

As noted in Section~\ref{s:stripping}, conditioning on $G\in\C(n,m)$ being simple is equivalent to taking $G\in\G(n,m)$.  For such $G$,
by Theorem~\ref{giantconc} and Lemma~\ref{kernelconc}, we have $\hat n \sim a_1n$ and $\hat m\sim a_2n$ a.a.s., for some constants $a_1$ and $a_2$ with $a_1<a_2$ and depending only on $c$. The same concentration then holds also for $G\in\C(n,m)$ by a quite simple argument: an alternative way to generate $G\in\C(n,m)$ is to first
decide how many loops, $\ell$, and multiple edges, $j$, it has (and their
multiplicities $m_1,\ldots, m_j$) with the correct probability,   generate an
underlying simple graph $G^*$ at random, and then adorn $G^*$ with $\ell$
loops at random locations, and the required extra copies of $j$ of its
edges. The distribution of $G^*$ should be uniform with $n$ vertices and
$m-\ell-\sum_i(m_i-1)$ edges, and the locations of the loops and multiple
edges are chosen at random.  Simple calculations with Markov's inequality show
that $\ell+\sum m_i=O(\log n)$ a.a.s.\ Furthermore, by
Lemma~\ref{treesizes}, adding a loop or giving an edge of $G^*$ extra
parallel copies, a.a.s.\ will not increase the size of its 2-core by more than
$O(\log n)$ (and of course cannot decrease it). Adding a loop or an extra copy
of an edge can only increase the kernel size by 2 vertices or 3 edges
(the extreme case is that of adding an edge parallel to an edge in the middle
of a path of vertices of degree 2). It follows that the concentration in Theorem~\ref{giantconc} and Lemma~\ref{kernelconc} also applies for $G\in\C(n,m)$, in particular $\hat n \sim a_1n$ and $\hat m\sim a_2n$ a.a.s.

{From} the above  paragraph, the numbers of vertices and edges of the 2-core of $G$ are a.a.s.\  $b(1-t)n+o(n)$ and $b_1n+o(n)$ respectively, and the number of degree 2 vertices is a.a.s.\ $b_2n+o(n)$. Recall that $b_1> b(1-t)>b_2$. Deleting isolated cycles to obtain $G_0$ as the trimmed core of $G$ will maintain uniform randomness, provided its numbers of vertices and edges are conditioned upon (and conditional upon having no isolated cycles). As observed by {\L}uczak~\cite{L}, the number of vertices in isolated cycles of the 2-core is small; it is easy to show that it is bounded in probability, or a.a.s.\ $O(\log n)$ for example. So these do not affect our argument and we ignore them.
The argument above shows that Lemma~\ref{treesizes} applies also to $G\in\C(n,m)$. It follows that the number of vertices of
$G_0$ that are adjacent to at least $N-1$ vertices of $G-G_0$ is a.a.s.\ at most $\eps_Nn$. This is therefore a bound on the number of red vertices in the painted kernel of $G_0$ with respect to $G$, and on the number of edges that are coloured red because they correspond to a maximal 2-path of $G_0$ containing such a vertex. Note that, given $G_0$ and the number $r$ of red vertices, each $r$-set of   vertices of $G_0$ is equally likely to be the set of red ones.

Similarly, since $b(1-t)>b_2$, the average number of vertices of $G_0$
assigned to an edge of its kernel is bounded, and,  using   Lemma~\ref{C32}, has an exponential tail. It follows that the number of edges of $K(G_0)$ that are painted red at the start of the stripping process is also at most $\eps_Nn$.

At the start of the stripping process, the kernel has $\hat n$ vertices. {From} Lemma~\ref{treesizes} and similar elementary analysis, the
proportion of its edges and vertices that are red is at most
$\eps_N$.  We need to consider how many red edges or vertices are produced
during the step in which
$G_{i+1}$ is obtained.  Purple edges and pink vertices merely denote ``potential trouble" and will
remain at the end, and we need their number to remain small. To aid in this, we define for fixed $\eps'>0$ the stopping time $T(\eps')$ to be the
 the smallest value $T$ such that at least one of the following holds:
\begin{itemize}
\item $T\ge \eps'\hat n$,
\item $K(G_T)$  has no red edges and no red vertices,
\item $K(G_T)$  has more than $\eps'n$ coloured edges and vertices.
\end{itemize}
Here $T$ is a stopping time for the exposition process, formally with respect
to the $\sigma$-fields generated by $R_i$, $P_i$ etc. We will examine the behaviour of the process up to the stopping time $T(\eps')$,  for $\eps'$ and $N$ fixed, and note the behaviour of the conclusions we draw, as $N\to\infty$.

We use Lemma~\ref{strip} frequently.
Note first that each stripping step deletes an absolutely bounded number of vertices and edges from the kernel. Hence by the definition of $T(\eps')$, for $i<T(\eps')$, there are at least $\hat n/2$ uncoloured edges in $K(G_i)$ (for $\eps'$ sufficiently small). Also, the number of new pink vertices or purple edges each increase  by at most 2 in each stripping step. So $K(G_i)$ has at most $4\eps'n$ such elements for $i<T(\eps')$.

We must also examine the distribution of numbers of red edges in $K(G_i)$ for $i<T(\eps')$. In Step $i+1$, first assume that the vertex
$u$ adjacent to the end of   the edge $\gamma_0$ to be removed  is
uncoloured and of degree 3.
The probability that either edge incident with $u$ was already
purple is, using Lemma~\ref{strip} and the bound on the number of purple edges, at most
$O(\eps')$ (where the implicit constant in this bound is independent of
$N$). So   $O(\eps')$ is an upper bound on the probability that a  new red
edge is created. On the other hand, the probability that $u$ is already
pink is $O(\eps')$  by Lemma~\ref{strip} and the conclusions in the above paragraph; then it becomes red and all incident
edges become red. The distribution of the number of such edges is
asymptotically truncated Poisson (the distribution is actually multinomial
conditioned on $|f^{-1}(j)|\ge 3$ for all $j$). This has an exponential
tail.

We conclude that there is an upper bound  $O(\eps')$ on the expected number  of new red edges arising in every step of the process, with an
exponential tail, up until time $T(\eps')$. This is, at each step, conditional upon the state of the process in the previous step. Note that the increase in the number of red edges in one step is $O(\log^2 n')$ with probability at least $1-o(n^{-2})$. A standard supermartingale inequality now shows that a.a.s.\ the total number of red edges created up to time $T(\eps')$ is $O(T(\eps')\eps' )$
(see~\cite[Corollary 4.1]{des}; the last paragraph of the proof of Theorem~5.1 in that reference explains how to handle the fact that the expected change in the number of red edges is not bounded by a constant). As $T(\eps')\le \eps'\hat
n+1$ by definition, there are a.a.s.\ at most   $O((\eps')^2\hat n)$ new red edges. Moreover, since every step of the process uses up a red edge (of which there are initially at most  $\eps_N \hat n$) this implies that $T(\eps')=O (\eps_N +(\eps')^2)\hat n$. Now   take a concrete instance  $g(N)$ of this function $\eps_N$, and   let  $\eps'=\eps'(N)=\sqrt{g(N)}$. Then $T(\eps')= o (\eps' \hat n)$ as $N\to\infty$. Hence, a.a.s.\ in the definition of $T(\eps')$ it is the condition that there are no red edges nor vertices that is the binding one.  That is, a.a.s.\ the whole process lasts for at most $\eps'(N)\hat n$ steps, where $\eps'(N)\to 0$ as $n\to\infty$.
Since $R_N(G)=G_M$, the lemma follows.
\qed
\smallskip

We denote by $G-E(H)$ the spanning subgraph of $G$ with edge set
$E(G)\setminus E(H)$. Note that no vertex of $R_N(G)$ has more than one
edge to $G_0- R_N(G)$, since otherwise all its  edges would be painted
red and so the vertex must be deleted eventually. Similarly,  when $G_0$ is
the 2-core of a graph $G$,  the edges incident with vertices of
$G_0$ that are adjacent to at least $N-1$ vertices  of $G-G_0$
are  initially made red. There are $O(\eps_N n)$ of these, and these
vertices cannot survive in $G_M$. It follows that no vertex of $R_N(G)$
has more than $N$ edges to
$G- R_N(G)$.

To prove Theorem~\ref{maincore} we also need to check the condition
on sizes of components. We first need a preparatory lemma on a sort of coalescing
branching process. This is a simplified version that does not apply directly to the actual process we need to consider, but will, with appropriate choice of $Z$, provide a useful comparison via stochastic domination.

\begin{lemma}\lab{coalesce}
Fix  $0<\eps<1/2$
and a nonnegative random variable $Z$ with $\ex Z<\eps$ and with exponential tail. Suppose that a graph $F_n$ with $n$ vertices has component sizes
with exponential tail. Suppose furthermore that at most $\delta n$ vertices
are marked, and that these include all the vertices in nontrivial
components. Now process the marked vertices successively by adding edges to
a random set of neighbours. The random number of neighbours chosen is distributed according to $Z$, and,  given the number
of neighbours, the neighbours themselves are chosen uniformly at random from all vertices. (These choices are done independently at random for each vertex processed. For simplicity, we permit loops, so a vertex may choose itself.)  Each vertex processed becomes unmarked, and each unmarked
isolated vertex that is joined to becomes marked. The process finishes
with a final graph, $F_n^{(f)}$
 when all vertices are unmarked. Then, for
$\delta$ sufficiently small, a.a.s.\ the component sizes in   $F_n^{(f)}$
have an exponential tail,  and a.a.s.\ there are at most $\eps n$ vertices in nontrivial components.
\end{lemma}

\noindent
  Here   $\delta$ only has to be smaller than some
absolute constant $c$.

\proof
Let us perform an equivalent process, in two stages: first growth, and secondly identification and pruning. Let $S$ denote the set of marked vertices. In the first stage, for each vertex $v\in S$, perform a Galton-Watson branching process with birth law given by the distribution of $Z$, and originating with the single individual $v$. All the children in all these  processes are at this point represented as separate vertices in a set $T$, where $S\cap T=\emptyset$, and the branching processes are represented as trees.

In the second stage, perform random identifications of vertices in the trees generated in the first stage: each vertex $u$ in $S$ is taken in turn, and for each such $u$, each child vertex $w$ in its branching process is taken in turn (working up the tree away from $u$). Then, with probability $k/n$, the vertex $w$ is identified with one of the vertices previously processed in this second stage, where $k$ is the number of such previous vertices. If identification occurs, the vertex to identify with is picked at random.
Furthermore,  all vertices
in the branches of the tree  above $w$  are deleted.

If we now add all edges present in $F_n$, it is clear that we obtain a graph with the same distribution as $F_n^{(f)}$.

To bound the component sizes of $F_n^{(f)}$, we analyse the equivalent process
{\em without} the deletion steps. Start with a vertex randomly chosen in
$S$. The size of its component in $F_n$ has exponential tail (in the
probabilistic sense) and has expected size    $1+O(\eps)$.    Each vertex in its
component is in $S$, and we may consider the  tree of each one separately. The
size of each of these trees has an exponential tail, and expected size
  $1+O(\eps)$. Each vertex in the tree is   identified with a number of
vertices processed earlier or later, and the number of these has an
exponential tail with expected size $O(\delta)$.  So we may consider a new
branching process, the children of a vertex being the new vertices in any tree
reached   by identification. The number of children then has an exponential
tail, and hence,  using~\cite[Theorem 3.3]{SW}, so does the size of the new
branching process. Furthermore, it is easy to see that its expected size is
  $1+O(\eps+\delta)$.

Because no truncations occur, the sizes of the new branching processes are independent and identically distributed. The rest of the proof is straightforward (c.f.\ Lemma~\ref{treesizes}).
\qed

To apply this lemma, the initially marked vertices are the initially red edges. New marked vertices are new red edges.

\begin{lemma} \lab{tails}
Let $G$ be as in Lemma~\ref{epsilons}. For $N$ sufficiently large, the set of numbers of
edges in  the components of
$G-E(R_N(G))$ a.a.s.\ has an exponential  tail.
\end{lemma}

\proof
Define $G_0$ as in  the proof of Lemma~\ref{epsilons}. It was shown there that
 the proportion of edges and vertices of the painted kernel  that are red  is at most $\eps_N$.

The number of new red edges generated in any one
step, as observed above, has an exponential tail with a truncated Poisson approximation.
We need to consider the components  $C_1,\ldots, C_j$ induced by all the
edges of  the kernel $K(G_0)$ that are red or are subsequently painted red.
We may begin with the initially red edges and vertices, all considered as marked vertices in some graph $F_q$ as in Lemma~\ref{coalesce}. These are arranged  in components of $F_q$ according to the  respective components in the subgraph of   $K(G_0)$ that they induce. These component sizes are easily seen to have an exponential tail. This follows because of two things: firstly, as noted in the proof of Lemma~\ref{epsilons}, the red vertices occur as a set chosen uniformly at random, and similarly the red edges that are not adjacent to red vertices. The second ingredient is that the distribution of degrees of the red vertices will be determined by the distribution of degrees of the vertices in $G_0$, which are, by the results in~\cite{CW}, multinomial and hence in the limit Poisson. We omit some details, as this part of the proof is straightforward.

  The other vertices and edges of $K(G_0)$ are all vertices of $F_q$.  We may assume as above that the number of red edges is always at most $\eps_Nn$. As each red edge is processed, its vertex in $F_q$ joins with  at most  probability $\eps_N$  to one or two other vertices of $F_q$. It is thus seen that the sizes of components in the resulting graph are bounded above by those of an associated process of the type analysed in Lemma~\ref{coalesce}, and the variable $Z$ has $\ex Z<\eps_N <1$.

{From} Lemma~\ref{coalesce}  we deduce
 that the sizes of $C_1,\ldots, C_j$ a.a.s.\ have an exponential tail. From Lemma~\ref{C32} it is easy to see that the set of numbers of degree 2 vertices in the maximal
2-paths of the 2-core of $G$ a.a.s.\ have an exponential tail (and this applies equally well if the isolated cycles are regarded as maximal 2-paths).  Finally,  Lemma~\ref{treesizes} says that the set of sizes trees attached to each vertex of the 2-core of $G$ has an exponential tail a.a.s. Recall that these are attached randomly.   Combining these statements using~\cite[Theorem 3.3]{SW} gives the result. \qed
  \smallskip

As we shall show formally at the end of  the next section, all that remains to prove Theorem~\ref{maincore} is to verify that the
$N$-reduced core a.a.s.\ satisfies the required   expansion property. This is considered in the next section.

\section{Expansion of the kernel} \lab{kernel expansion}

 Let $\G(\dv)$ denote the uniform  probability space of the graphs with
degree sequence $\dv = (d_1,\ldots, d_{n})$. To model this probability
space we may use the pairing model (see~\cite{models}). Here there are
cells labelled
$1,\ldots, n$ with $d_i$ points in the $i$th cell. A uniformly random
pairing of all the points is selected, denoted $\P(\dv)$. Regarding the
cells as vertices, this produces a random (pseudo-)graph $G$ on $n$
vertices. It is easy to check that the following holds (see~\cite{CW} for example).
\begin{lemma}\lab{equiv}
  The distribution of   $G$ arising from $\P(\dv)$ is exactly the same as that obtained by  restricting
$\C(n,m)$ to graphs with degree sequence $\dv$.
 \end{lemma}

 To make the connexion with the   results in Section~\ref{s:stripping} we still  need the following.
\begin{lemma}\lab{simple} For $m=O(n)$, there exists $C>0$ such that
  $G\in \C(n,m)$ is simple with probability at least $C$.
\end{lemma}
\proof It is straightforward to show using the method of moments that the
numbers of loops and pairs of parallel edges are asymptotically independent
Poisson in distribution. The expected numbers are $m/n$ and ${m \choose
2}/n^2$ respectively, which are both $O(1)$. The result follows.
\qed

So we may focus on the pairing model, and Lemma~\ref{epsilons} tells us that we only need to consider a
set of degree sequences that a.a.s.\ contains the degree sequence of
$\C_3(n,m)$, for $m=O(n)$. The convergence expressed in the following lemma is uniform over all degree sequences $\dv$ in the stated range.

\begin{lemma}\lab{kernel-expander}
 For   some $\alpha>0$, the random multigraph   $G_\dv$ arising from the
pairing model $\P(\dv)$ with $3\le \min d_i\le \max d_i \le n^{0.02}$   is
a.a.s.\ an
$\alpha$-expander.
\end{lemma}
\proof First, note that
$\Phi$ is the minimum of $e(S)/d(S)$  over all sets $S$ of vertices whose
sum of degrees is at most the number of edges  of the graph, where $d(S)$ is
the sum of degrees of vertices in $S$ and $e(S)$ is the number of edges
leading out of $S$.

In the proof of~\cite[Lemma 12.6]{L}, {\L}uczak shows  that, conditional
on  a given degree sequence $\dv$ with minimum   3 and maximum   at most
$n^{0.02}$, the multigraph arising in the random pairing $\P(\dv)$  a.a.s.\
has no subgraph on
$r$ vertices,
$2\le r\le n^{0.4}$, with more than $1.2r$ edges. It follows that each set
$S$ of at most $n^{0.4}$ vertices has $e(S)/d(S)\ge 1/5$. Much simpler
calculations show that the same is true for $r=1$:   the expected number of
pairs of loops with the same vertex is $o(1)$.

For the sets of vertices between $n^{0.4}$ and $n/2$, {\L}uczak only
establishes a constant lower bound on the number of edges leaving the set.
Here we more than fill the gap by computing the expected number of sets of
vertices $S$ with $d(S)=q$, where $n^{0.2}<q\le m$ ($m=|E(G_3)|$) and
$e(S)=t$. Note that the upper bound $\max d_i\leq
n^{0.02}$ is not required in this part.

Assume $|S|=s$; we will sum over all relevant $s$ later. The  expected
number of   sets as above is

\bel{prob}
\sum_* P(m,t,q)
\le {n\choose s}P(m,t,q)
\ee
 where
 the summation is over all subsets $S$ of $V$ with $|S|=s$ and
$d(S)=q$, and $P(m,t,q)$ is the probability that a random matching of
$2m$ points has exactly
$t$ edges leaving a given set of $q$ points. Thus
$$P(m,t,q) = {2m-q\choose t}{q\choose t} t!M(q-t)M(2m-q-t)M(2m)^{-1} $$
where  $M(k)=\frac{k!}{(k/2)!2^{k/2}}$ is the number of perfect matchings of $k$
points ($k$ even).

Using Stirling's formula and separating out insignificant factors of size
$m^{O(1)}$ we obtain
 \bean P(m,t,q) &\le& m^{O(1)} \frac  {(1-\frac{q}{2m})^{2m-q}(\frac
{q}{t})^t(\frac {q}{2m})^{q-t}   } { (\frac {q-t}{2m})^{(q-t)/2} (1-\frac
{q+t}{2m})^{m-(q+t)/2} }\\ &\le &
m^{O(1)}(O(q/t))^t\left(1-\frac{q}{2m}\right)^{ m-q/2}
\left( \frac{q}{2m}\right)^{q/2-t/2}\\ &\le &  m^{O(1)}\left(
(O(1/\eps))^\eps\left(r^{1-\eps}(1-r)\right)^{ 1/2}
 \right)^q
\eean using $q\le m$ and with $t=\eps q$ and $r=q/(2m)$. By taking $\eps$
close to 0 we can make $\eps^\eps$ close to 1, and so the most significant
part of this   is
\bel{part2} (r^{1/2}(1-r)^{1/2})^q.
\ee

First consider the case that  $q/3\le n/2$.   Since $s\le q/3$ by the fact
that vertex degrees are all at least 3, we can use
$${n\choose s}\le {n\choose  q/3} \le (3en/ q)^{q/3}\le  (2em/ q)^{q/3} = (e/ r)^{q/3}. $$  Multiplying by~\eqn{part2}, we maximise
$r^{1/6}(1-r)^{1/2}$ at $r=1/4$ and find that~\eqn{prob} is at most
\bel{part3} m^{O(1)} ( 1-\eps')^q
\ee for some $\eps'>0$ when $\eps <\eps_0$ (some $\eps_0$  sufficiently
small). On the other hand, if $q/3> n/2$, use  ${n\choose s}\le 2^n\le
2^{2q/3}$ and since $r^{1/2}(1-r)^{1/2}\le 1/2$ the same conclusion is
reached. Then summing~\eqn{part3} over all $n^{0.2}<q\le m$, all relevant
$s$ and  all
$t\le
\eps_0q$, the result is $o(1)$. Hence the expected number of sets of
vertices in the size range being considered ($n^{0.2}\le s \le n/2$) with
$e(S)/d(S)=t/q<\eps_0$ is o(1).

We conclude that  a.a.s.\  $\Phi\ge \min(1/5,\eps_0)$.
\qed
\smallskip

\noindent
{\bf Proof of Theorem~\ref{maincore}\ } For $m$ as in the theorem statement, let    $G\in\C(n,m)$. Then, by Lemma~\ref{epsilons}, for some constant $c'>0$ the number $n'$  of vertices in $K(R_N(G))$  is a.a.s.\ at least $c'n$  for $N$ sufficiently large. Moreover this graph is distributed as    $\C_3(n',m')$, given $n'$ and its number of vertices $m'$. Then by Lemma~\ref{equiv}, further restricting this to degree sequence $\P(\dv)$ gives graphs with the distribution of $\P(\dv)$. It is well known and easy to verify  that, since $m'=O(n')$, a.a.s.\ the maximum degree occurring in $\C_3(n',m')$ is $o(n^{0.02})$. Hence by Lemma~\ref{kernel-expander}, $K(R_N(G))$ is a.a.s.\ an  $\alpha'$-expander (with $\alpha'$ being the $\alpha$ from   that lemma). Thus $R_N(G)$ is a.a.s.\ an $\alpha $-expander, where $\alpha=\alpha'/N$ say, since it is obtained from its kernel by inserting at most $N$ vertices of degree 2 into each edge.   It a.a.s.\ satisfies property (2) in the definition o!
 f an $\alpha $-strong core of $G$ by Lemma~\ref{tails}. It satisfies property (3) for $\alpha<1/(2N)$  by the definition of the severe stripping process, since vertices of the 2-core that are adjacent to more than $N$ vertices outside it are deleted, and during the stripping, any vertices adjacent to at least two that have been deleted during stripping are deleted themselves.  Thus for $\alpha$ sufficiently small, $R_N(G)$ is a.a.s.\ an $\alpha $-strong core of $G$. The theorem then follows by Lemma~\ref{simple}, which lets us translate results holding a.a.s.\ to $\G(n,m)$. \qed

 \medskip

 \noindent
 {\bf Acknowledgment} The authors wish to thank Elchanan Mossel for helpful discussions at an early stage of this research.

\end{document}